\documentclass[12pt]{amsart}
\usepackage{amsfonts}
\usepackage{amssymb}
\usepackage{graphicx}
\def\({\left(}       \def\){\right)}

\newtheorem{prop}{\sc Proposition}
\newtheorem{lem}{\sc Lemma}
\newtheorem{theorem}{\sc Theorem}
\newtheorem{cor}{\sc Corollary}

\newenvironment{pf}{\noindent{\textit{Proof. }}}{$\Box$ }

\begin{document}
%%%%%%%%%%%%%%%%%%	Title, etc.	 %%%%%%%%%%%%%%%%%%
\title[On the harmonic M\"obius transformations]{On the harmonic M\"obius transformations}
\author[R. Hern\'andez]{Rodrigo Hern\'andez}
\address{Facultad de Ingenier\'{\i}a y Ciencias, Universidad Adolfo Ib\'a\~nez, Av. Padre Hurtado 750, Vi\~na del Mar, Chile.} \email{rodrigo.hernandez@uai.cl}

\author[M. J. Mart\'{\i}n]{Mar\'{\i}a J. Mart\'{\i}n}
\address{Department of Physics and Mathematics, University of Eastern Finland, P.O. Box 111, FI-80101 Joensuu, Finland.} \email{maria.martin@uef.fi}

\thanks{This research is supported by grant Fondecyt $\#1150284$, Chile. The second author also thankfully acknowledges partial support from Academy of Finland grant $\#286877$, Spanish MINECO/FEDER research project MTM2015-65792-P, and Spanish Thematic Research Network MTM2015-69323-REDT, MINECO}

\begin{abstract} It is well-known that two locally univalent analytic functions $\varphi$ and $\psi$ have equal Schwarzian derivative if and only if there exists a non-constant M\"obius transformation $T$ such that $\varphi=T\circ \psi$. In this paper, we identify completely the relationship between two locally univalent harmonic mappings with equal (harmonic) Schwarzian derivative.
\end{abstract}

%%%%%%%%%%%%%%%%%%%%%%%%%%%%%%%%%%%%%%%%%%%%%%%%%%%
\date{\today}
\maketitle

\section*{Introduction}

The \emph{Schwarzian derivative} of a locally univalent analytic function $\varphi$ in a simply connected domain $\Omega$ of the complex plane $\mathbb{C}$ equals
\begin{equation}\label{eq-classicalSchwarzian}
S(\varphi)=\left(\frac{\varphi''}{\varphi'}\right)' -\frac 12\left(\frac{\varphi''}{\varphi'}\right)^2\,.
\end{equation}
The quotient $\varphi''/\varphi'$, denoted by $P(\varphi)$, is known as the \emph{pre-Schwarzian derivative} of the function $\varphi$.
\par\smallskip

A complex-valued function $f$ is said to be \emph{harmonic} in a simply connected domain $\Omega\subset \mathbb{C}$ if both $Re f$ and $Im f$ are real harmonic in $\Omega$. Every such $f$ can be written in the form $f=\overline g+h$, where both $g$ and $h$ are analytic in $\Omega$  (see \cite[p. 7]{Dur-Harm}). \par
It is known that the (second complex) \emph{dilatation} $\omega=g'/h'$ of a harmonic function $f=\overline{g}+h$ stores important information about $f$. For example, a non-constant harmonic mapping is \emph{orientation-preserving} if and only if $|\omega|\leq 1$ \cite[p. 8]{Dur-Harm}. Lewy \cite{Lewy} proved that a necessary and sufficient condition for $f$ to be locally univalent is that its \emph{Jacobian} $J_f=|f_z|^2-|f_{\overline z}|^2=|h'|^2-|g'|^2$ is different from $0$. Hence, a locally univalent harmonic mapping $f$ is orientation-preserving if its Jacobian is positive and \emph{orientation-reversing} if $J_f<0$. Note that the locally univalent harmonic mapping $f=\overline{g}+h$ is orientation-preserving if and only if $h$ is locally univalent and the dilatation $\omega=g'/h'$ is an analytic function bounded (in modulus) by $1$.
\par\smallskip
The \emph{harmonic Schwarzian derivative} $S_H$ of a locally univalent harmonic function $f$ with Jacobian $J_f$ was defined in \cite{HM-Schwarzian} by
\begin{equation}\label{eq-Schwarzian0}
S_H(f)=\frac{\partial}{\partial z}\left(P_H(f)\right)-\frac 12 \left(P_H(f)\right)^2\,,
\end{equation}
where $P_H(f)$ is the \emph{harmonic pre-Schwarzian derivative} of $f$, which equals
\[
P_H(f)=\frac{\partial}{\partial z} \log J_f\,.
\]
It is easy to check that $S_H(f)=S_H(\overline f)$ for any locally univalent harmonic mapping $f$ in a simply connected domain $\Omega$. Therefore, without loss of generality we can assume that $f$ is orientation-preserving.  The harmonic Schwarzian derivative of the orientation-preserving harmonic mapping $f=\overline g+h$ with dilatation $\omega=g^\prime/h^\prime$ can be written as
\begin{equation}\label{eq-Schwarzian}
S_H(f)=S(h)+\frac{\overline \omega}{1-|\omega|^2}\left(\frac{h''}{h'}\,\omega'-\omega''\right)
-\frac 32\left(\frac{\omega'\,\overline
\omega}{1-|\omega|^2}\right)^2\,,
\end{equation}
where $S(h)$ is the classical Schwarzian derivative of the function $h$ defined by \eqref{eq-classicalSchwarzian}.
\par\smallskip
It is clear that if $f$ is analytic (so that its dilatation is identically zero) then $S_H(f)$ coincides with the classical Schwarzian derivative of $f$. In other words, the operator defined by \eqref{eq-Schwarzian0} -or by \eqref{eq-Schwarzian}- can be considered as a generalization of the classical Schwarzian derivative~\eqref{eq-classicalSchwarzian}. We refer the reader to \cite{HM-Schwarzian} for the motivation of the definition -related to the classical argument of approximation by M\"{o}bius transformations that seems to go back to E. Cartan \cite{Cartan} (see also \cite[p. 113]{Gardiner})- as well as for different properties that the harmonic Schwarzian derivative satisfies. Some of them -the most important ones for our purposes- will be considered in Section~\ref{ssec-properties} below.
\par\smallskip
The harmonic Schwarzian operators $P_H$ and $S_H$ defined above have proved to be useful to generalize classical results in the setting of analytic functions to the more general setting of harmonic mappings. For instance:
\par
- The classical criteria of univalence and quasi-conformal extension for analytic functions in terms of the pre-Schwarzian derivative due to Becker \cite{Becker} (see also \cite{Becker-Pom-1}) and Ahlfors \cite{Ahlfors} are generalized to those cases when the functions considered are merely harmonic in \cite{HM-qc} and \cite{HM-Schwarzian}.
\par
- In the article \cite{HM-Nehari}, the celebrated criterion of univalence in terms of the Schwarzian derivative obtained by Nehari \cite{Nehari} as well as the corresponding criterion for quasi-conformal extension due to Ahlfors and Weill \cite{AW} are generalized to the harmonic setting.
\par
- Two criteria for the bounded valence of harmonic mappings in terms of the harmonic pre-Schwarzian and Schwarzian derivatives, respectively, are obtained in \cite{Huusko-M} as a generalization of some of the results in \cite{Becker-Pom-2} and \cite{GP}.
\par
- The relationship between John disks and the pre-Schwarzian derivative analyzed in \cite{Hag-Hag} (see also \cite{Ch-O-P}) is extended to the harmonic setting in \cite{Chen-Ponn}.
\par\smallskip
These harmonic operators have recently attracted the attention of several authors. We refer the reader to the recent papers \cite{Graf, Liu-Ponn, Liu-Ponn-2}, for instance. In fact, also recently, it has been proved that the harmonic Schwarzian derivative defined above would help to solve some related problems regarding harmonic mappings in the unit disk. More concretely, a consequence of the main theorem in \cite{CHM} is that if it would be possible to show that the norm of the harmonic Schwarzian derivative of any univalent harmonic mapping in the unit disk $\mathbb D$ is bounded by $19/2$, then the order of the well-known family of suitable normalized orientation-preserving univalent harmonic functions in $\mathbb D$ would be equal to 3, as conjectured (see \cite{C-SS}, \cite{SS}, or \cite[Sec. 5.4]{Dur-Harm}).
\par\smallskip
However, there are still fundamental questions related to the harmonic Schwarzian derivative that remain unresolved. For instance, as mentioned in Nehari's paper \cite{Nehari}, a locally univalent analytic function $\varphi$ with Schwarzian derivative $S(\varphi)$ necessarily equals the quotient $u_1/u_2$, where $u_1$ and $u_2$ are linearly independent solutions of the linear differential equation
\[
u''+\frac{S(\varphi)}{2}u=0\,.
\]
\par
This readily implies that if two given functions $\varphi_1$ and $\varphi_2$ (again locally univalent and analytic) have equal Schwarzian derivative, then $\varphi_1=T\circ \varphi_2$, where $T$ is a non-constant \emph{M\"obius transformation} $T$ (also called \emph{linear fractional transformation}) of the form
\begin{equation}\label{eq-Mobius}
T(z)=\frac{az+b}{cz+d}\,,\quad z \in \mathbb C\quad \text{and}\quad ad-bc\neq 0\,.
\end{equation}

A straightforward calculation shows that if $\varphi_1=T\circ \varphi_2$ then $S(\varphi_1)=S(\varphi_2)$. Hence we obtain that  $\varphi_1$ and $\varphi_2$ have equal Schwarzian derivative if and only if $\varphi_1=T\circ \varphi_2$ for some non-constant M\"obius transformation $T$.
\par\smallskip
The main purpose of this paper is to characterize the relationship between two locally univalent harmonic mappings with equal harmonic Schwarzian derivative.
\par\smallskip
It sometimes happens that in order to get the solution of certain problem, it is needed to distinguish between those cases when the functions $f=\overline{g}+h$ considered have constant dilatation $\omega=g'/h'$ and those cases when the dilatations are not constant (see, for instance, \cite{CM}). This is precisely what occurs in our case, so that we need to consider these situations separately. Notice that the dilatation of a locally univalent harmonic mapping is not constant if and only if the functions $g'$ and $h'$ in the decomposition $f=\overline g+h$ are linearly independent.
\par\smallskip
The paper is organized as follows. In Section~\ref{sec-background}, we review some of the properties of the operator $S_H$ and use them to normalize the functions considered in a suitable way. In Section~\ref{sec-constantdilatation}, we completely identify the relationship between two harmonic mappings with equal harmonic Schwarzian derivatives in those cases when the dilatation of one of the functions involved is constant. In Section~\ref{sec-main}, we treat the more difficult problem of determining the relation between two harmonic maps $f_1$ and $f_2$ with $S_H(f_1)=S_H(f_2)$ when the dilatations of both $f_1$ and $f_2$ are non-constant. Finally, in Section~\ref{sec-solution}, we state the complete solution to the problem considered.

\section{Background}\label{sec-background}

We would like to start this section with some comments related to the Schwarzian derivative operator $\mathcal S$ introduced by Chuaqui, Duren, and Osgood in \cite{Ch-D-O} for the family of harmonic functions $f=\overline g+h$ such that $\lambda_f=|h'|+|g'|\neq 0$ and whose dilatation $\omega=q^2$ for some meromorphic function $q$. The precise definition of $\mathcal S(f)$ for such harmonic mappings $f$ is
\[
\mathcal S(f)=2\left(\frac{\partial^2}{\partial z^2} \lambda_f -\left(\frac{\partial}{\partial z} \lambda_f\right)^2\right)\,.
\]

It is well-known that those harmonic mappings that satisfy the specified conditions (which are not necessarily locally univalent) can be lifted to a minimal surface \cite[Ch. 10.2]{Dur-Harm}. By exploiting this connection with differential
geometry, Chuaqui, Duren, and Osgood have obtained many interesting results on different properties for the \emph{lifts} of the given harmonic mappings in terms of $\mathcal S$: criteria for univalence (see \cite{Ch-D-O-univalence}), quasi-conformal extension \cite{Ch-D-O-qcextension}, or distortion theorems (see \cite{Ch-D-O-distortion} and also \cite{Ch-D-O-M-M-M-O}). In fact, in \cite{Ch-D-O}, using again techniques from differential geometry, the authors are able to characterize the relationship between two harmonic functions $f_1$ and $f_2$ that satisfy the specified conditions and such that $\mathcal  S(f_1)=\mathcal  S(f_2)$.
\par\smallskip
It is not clear to us how to apply similar differential geometry tools to the solution of the problem that we consider in this article, since the dilatations of locally univalent harmonic functions $f=\overline{g}+h$ do not necessarily equal the square of a meromorphic function. Therefore, a different approach has been developed in order to determine the relationship between two locally univalent harmonic mappings with equal harmonic Schwarzian derivative $S_H$ as in \eqref{eq-Schwarzian}. In fact, perhaps it is convenient to state that our approach cannot be used to provide an alternative proof to the characterization of the relationship between two harmonic functions $f_1$ and $f_2$ that satisfy the specified conditions above and such that $\mathcal  S(f_1)=\mathcal  S(f_2)$, obtained in \cite{Ch-D-O}. The reason is that one of the main tools we will use is the invariance of $S_H$ under pre-composition with locally univalent \emph{affine} harmonic functions, explained in Section~\ref{ssec-properties} below, which is a property that is not shared by the harmonic Schwarzian operator $\mathcal S$.
\par\smallskip

\subsection{Some properties of the harmonic Schwarzian derivative}\label{ssec-properties}

It was shown in \cite{HM-Schwarzian} that for any given orientation-preserving harmonic mapping $f=\overline{g}+h$ with dilatation $\omega$ in a simply connected domain $\Omega$ and all $z_0\in\Omega$
\begin{equation}\label{eq-Schwarzian2}
S_H(f)(z_0)=S(h-\overline{\omega(z_0)}g)(z_0)\,.
\end{equation}
\par
A straightforward calculation shows that whenever $f$ is a locally univalent harmonic mapping and $\phi$ is an analytic function such that the composition $f\circ\phi$ is well-defined, the harmonic Schwarzian derivative of the harmonic mapping $f\circ\phi$ can be computed using the \emph{chain rule}:
\[
S_H(f\circ\phi)=(S_H(f)\circ\phi)\cdot(\phi')^2+S(\phi)\,.
\]
\par
Another important property verified by the harmonic Schwarzian derivative $S_H$ is the invariance with respect to pre-compositions with affine harmonic mappings. More specifically, consider a \emph{locally univalent affine harmonic mapping}
\[
A(z)=a\overline z+b z +c\,,
\]
where $a$, $b$, and $c$ are complex numbers with $|a|\neq |b|$. Then for any locally univalent harmonic mapping $f$ the composition $A\circ f$ is harmonic and locally univalent as well and $S_H(A\circ f)=S_H(f)$. In particular (setting $a=1$ and $b=c=0$), we have $S_H(\overline f)=S_H(f)$. That is, as was mentioned before, we can assume without loss of generality that the functions considered are orientation preserving.
\par
Let $f=\overline g+h$ be an orientation preserving mapping and let $\mu$ have modulus one.  Define the \emph{anti-analytic rotation}
\begin{equation}\label{eq-antianalyticrotation}
R_\mu(f)=\mu\overline g+h\,.
\end{equation}
Notice that if $\omega$ denotes the dilatation of $f$, then the dilatation $\omega_\mu$ or $R_\mu(f)$ is $\omega_\mu=\overline\mu \omega$. It is obvious from \eqref{eq-Schwarzian0} that $S_H(R_\mu\circ f)=S_H(f)$.
\par\smallskip
The following result was proved in \cite{HM-Schwarzian}.

\begin{theorem}[\cite{HM-Schwarzian}, Thms. 1 and 2]\label{thm-constantdilatation} Let $f$ be a locally univalent (orientation-preserving, without loss of generality) harmonic function in a simply connected domain $\Omega$. Then, $S_H(f)$ is harmonic if and only if the dilatation of $f$ is constant. That is, if and only if $f=\alpha\overline h+h+\gamma$ for some constants $\alpha$ and $\gamma$ with $|\alpha|\neq 1$ $(|\alpha|<1)$ and some locally univalent analytic function $h$ in $\Omega$. In this case, $S_H(f)=S(h)$ so that $S_H(f)$ is in fact analytic.
\end{theorem}

\subsection{A useful normalization}\label{ssec-normalization} The properties of the Schwarzian derivative operator $S_H$ stated in the previous section allow us to make certain normalizations that will be useful to determine the relationship between two harmonic functions $f_1=\overline{g_1}+h_1$ and $f_2=\overline{g_2}+h_2$ on a simply connected domain $\Omega$ that satisfy $S_H(f_1)=S_H(f_2)$.
\par
Recall that we can assume that both harmonic mappings $f_1$ and $f_2$ preserve the orientation. Assume that $S_H(f_1)=S_H(f_2)$ is not harmonic so that, according to Theorem~\ref{thm-constantdilatation} above, the dilatations $\omega_1=g_1'/h'_1$ and $\omega_2=g'_2/h'_2$ (which turn out to be analytic functions in $\Omega$) are not identically constant.
\par
Consider any conformal mapping $\phi$ from the unit disk $\mathbb{D}$ onto $\Omega$ with $\phi(0)=z_0\in\Omega$ and use the chain rule for the harmonic Schwarzian derivative, as well as the invariance under pre-composition with affine harmonic mappings, to get that the functions $f_1\circ\phi-f_1(z_0)$ and $f_2\circ \phi-f_2(z_0)$, now sense-preserving harmonic mappings in the unit disk that fix the origin, have equal harmonic Schwarzian derivatives and state the problem in the unit disk. An easy calculation shows that the dilatations or these new harmonic functions are, respectively, $\omega_1\circ\phi$ and $\omega_2\circ\phi$. In order not to burden the notation, we again use $f_1=\overline{g_1}+h_1$ and $f_2=\overline{g_2}+h_2$ to denote the new mappings and $\omega_1$ and $\omega_2$ to denote the corresponding dilatations, that are non-constant analytic functions in $\mathbb{D}$.
\par
Since multiplication by a non-zero constant does not change the Schwarzian derivative, we can also suppose that $h'_1(0)=h'_2(0)=1$.
\par
Bearing in mind that none of the dilatations $\omega_1$ or $\omega_2$ are constant, we have that there exists a point $w\in\mathbb{D}$ such that both $\omega'_1(w)$ and $\omega'_2(w)$ are different from zero. Using again the chain rule and the transformation
\begin{equation}\label{eq-aff1}
f=\overline g+h \mapsto\frac{f\circ\varphi_w-f(w)}{h'(w)(1-|w|^2)}\,,
\end{equation}
where $\varphi_w$ is the automorphism of the unit disk given by
\begin{equation}\label{eq-automorphism}
\varphi_w(z)=\frac{w+z}{1+\overline w z}\,,
\end{equation}
we see that without loss of generality $w$ can be taken to be the origin.
\par
Finally, as it is explained on \cite[Sec. 5.1]{Dur-Harm}, we can apply \emph{invertible} affine transformations (an operation that does not change the Schwarzian derivative) to both $f_1$ and $f_2$ to get new functions with dilatations that fix the origin and satisfy $\omega'_1(0)\neq 0$ and  $\omega'_2(0)\neq 0$.
\par
Summing up, these arguments show that in order to characterize the relationship between two orientation-preserving harmonic mappings $f_1$ and $f_2$ in a simply connected domain $\Omega$ having equal harmonic Schwarzian derivatives which is not a harmonic mapping, we can assume without loss of generality that $\Omega=\mathbb{D}$ and that the two harmonic mappings $f_1=\overline{g_1}+h_1$ and $f_2=\overline{g_2}+h_2$ in the unit disk with dilatations $\omega_1$ and $\omega_2$, respectively, satisfy
\begin{equation}\label{eq-normalizationh}
h_1(0)=g_1(0)=1-h'_1(0)=h_2(0)=g_2(0)=1-h'_2(0)=0\,,
\end{equation}
\begin{equation}\label{eq-normalizationw}
 \omega_1(0)=\omega_2(0)=0\,,
\end{equation}
$\omega'_1(0)\neq 0$, and $\omega'_2(0)\neq 0$. Indeed, by choosing appropriates $\mu_1$ and $\mu_2$ of modulus one, we can consider the anti-analytic rotations $R_{\mu_1}$ and $R_{\mu_2}$ defined as in \eqref{eq-antianalyticrotation}, and apply them to $f_1$ and $f_2$, respectively, to have that the value of the derivative at the origin of the dilatations is not only different from zero but a positive real number. In other words, we can suppose
\begin{equation}\label{eq-normalizationw'}
\omega'_1(0)>0\quad\text{and}\quad \omega'_2(0)> 0\,.
\end{equation}
\par
A careful reader may have notice that in order to get the normalizations mentioned, we assumed that the harmonic Schwarzian derivatives of both the harmonic mappings $f_1$ and $f_2$ are not harmonic. This is equivalent, according to Theorem~\ref{thm-constantdilatation}, to the fact that the (second complex) dilatations of the functions involved are not constant. Obviously, the harmonic mapping $f=\overline g + h$ has constant dilatation if and only if $g'$ and $h'$ are linearly dependent.
\par
As pointed out in the introduction, in order to solve the problem we consider, we will need to distinguish between those cases when the functions involved have constant dilatation (so that the harmonic Schwarzian derivatives are harmonic) and those cases when the dilatations are not constant. The easier case when $S_H(f_1)=S_H(f_2)$ is a harmonic mapping will be treated in the next section.

\section{The linearly dependent case}\label{sec-constantdilatation}

We now solve the problem of characterizing the relationship between two locally univalent harmonic mappings $f_1$ and $f_2$ for which $S_H(f_1)=S_H(f_2)$ is a harmonic function. Recall that we can assume, without loss of generality that both functions are orientation-preserving after taking complex conjugates if needed. Also, that those harmonic mappings for which their harmonic Schwarzian derivative is a harmonic function are completely described in Theorem~\ref{thm-constantdilatation}.
\par
\begin{theorem}\label{thm-constdilat1}
Let $f_1$ and $f_2$ be two orientation-preserving harmonic mappings in a simply connected domain $\Omega\subset\mathbb{C}$. Suppose that $S_H(f_1)=S_H(f_2)$ is a harmonic function, so that $f_1=\alpha_1\overline{h_1}+h_1+\gamma_1$, where $h_1$ is an analytic function in $\Omega$, $\alpha_1\in\mathbb{D}$, and $\gamma_1\in\mathbb{C}$. Then there exist a M\"{o}bius transformation $T$ as in \eqref{eq-Mobius} and two complex numbers $\alpha_2\in\mathbb{D}$ and $\gamma_2\in\mathbb{C}$ such that
\[
f_2=\alpha_2\ \overline{T\circ h_1}+T\circ h_1+\gamma_2\,.
\]
\end{theorem}
\begin{pf}
The identity $S_H(f_1)=S_H(f_2)$, the fact that these functions are harmonic, and Theorem~\ref{thm-constantdilatation} give $f_2=\alpha_2\overline{h_2}+h_2+\gamma_2$ for certain $\alpha_2\in\mathbb{D}$ and $\gamma_2\in\mathbb{C}$. Moreover, $S_H(f_1)=S(h_1)$ and $S_H(f_2)=S(h_2)$. Therefore, there exists $T$ as in \eqref{eq-Mobius} such that $h_2=T\circ h_1$. This proves the theorem.
\end{pf}
\par\smallskip
The remaining part of this paper will be devoted to solve the more difficult problem of determining the relation between two harmonic mappings $f_1$ and $f_2$ for which $S_H(f_1)=S_H(f_2)$ is \emph{not} a harmonic function.

\section{The case of linear independence}\label{sec-main}
In order to solve the problem of determining the relationship between two locally univalent harmonic mappings $f_1=\overline{g_1}+h_1$ and $f_2=\overline{g_2}+h_2$ in a simply connected domain $\Omega\subset\mathbb{C}$ with non-constant dilatations $\omega_1=g'_1/h'_1$ and $\omega_2=g'_2/h'_2$, respectively (so that $S_H(f_1)$ and $S_H(f_2)$ are not harmonic), and such that $S_H(f_1)=S_H(f_2)$ we can argue as in Section~\ref{ssec-normalization} to obtain that after applying a series of affine transformations, anti-analytic rotations, and a composition with a Riemann map from the unit disk onto $\Omega$, we can assume without loss of generality that $\Omega=\mathbb{D}$ and the normalizations \eqref{eq-normalizationh}, \eqref{eq-normalizationw}, and \eqref{eq-normalizationw'} are satisfied by both $f_1$ and $f_2$.

\subsection{Towards the solution of the problem}\label{ssec-equations}
We now prove one of the key results in this paper.
\begin{prop}\label{prop-schw h}
Let $f_1=\overline{g_1}+h_1$ and $f_2=\overline{g_2}+h_2$ be orientation-preserving harmonic mappings in the unit disk with non-constant dilatations $\omega_1$ and $\omega_2$, respectively. Assume that both $f_1$ and $f_2$ are normalized as in \eqref{eq-normalizationh}, \eqref{eq-normalizationw}, and \eqref{eq-normalizationw'}. Suppose that $S_H(f_1)=S_H(f_2)$. Then the following assertions hold.
\begin{itemize}
\item[i)] $S(h_1)=S(h_2)$.
\item[ii)] For all $z\in\mathbb{D}$,
\begin{equation}\label{eq-prop-aux1}
\omega'_1(0) \left(\frac{h''_1(z)}{h'_1(z)}\, \omega'_1(z)-\omega''_1(z)\right)=\omega'_2(0) \left(\frac{h''_2(z)}{h'_2(z)}\, \omega'_2(z)-\omega''_2(z)\right)\,.
\end{equation}
\item[iii)] The identity
\begin{align}\label{eq-prop-aux2}
\nonumber &\overline{\left(\frac{h''_1(0)}{h'_1(0)}\, \omega'_1(0)  -\omega''_1(0)\right)} \omega_1(z)-\frac 32 (\omega'_1(0))^2\, (\omega_1(z))^2\\
& \quad \quad \quad =\overline{\left(\frac{h''_2(0)}{h'_2(0)}\, \omega'_2(0)-\omega''_2(0)\right)} \omega_2(z)-\frac 32 (\omega'_2(0))^2\, (\omega_2(z))^2
\end{align}
holds for all $z$ in the unit disk.
\end{itemize}
\end{prop}
\begin{pf}
Since the functions $h_1$ and $h_2$ are locally univalent mappings in the unit disk, the (classical) Schwarzian derivatives $S(h_1)$ and $S(h_2)$ of these functions are analytic in $\mathbb{D}$. Hence, to show that $S(h_1)=S(h_2)$ it suffices to see that the Taylor coefficients of $S(h_1)$ and $S(h_2)$ are equal.
\par
The condition that $S_H(f_1)=S_H(f_2)$ is, by \eqref{eq-Schwarzian}, equivalent to
\begin{align}\label{eq-prop-equalSh}
\nonumber S(h_1)& +\frac{\overline{\omega_1}}{1-|\omega_1|^2}\left(\frac{h''_1}{h'_1}\, \omega'_1-\omega''_1\right)-\frac 32 \left(\frac{\omega'_1\overline{\omega_1}}{1-|\omega_1|^2}\right)^2\\
&=S(h_2) +\frac{\overline{\omega_2}}{1-|\omega_2|^2}\left(\frac{h''_2}{h'_2}\, \omega'_2-\omega''_2\right)-\frac 32 \left(\frac{\omega'_2\overline{\omega_2}}{1-|\omega_2|^2}\right)^2\,,
\end{align}
which implies that $S(h_1)(0)=S(h_2)(0)$ since  the dilatations $\omega_1$ and $\omega_2$ of $f_1$ and $f_2$, respectively, fix the origin.
\par
A straightforward calculation shows
\begin{align*}
\frac{\partial S_H(f_1)}{\partial z}&=(S(h_1))'+\overline{\omega_1}\frac{\partial}{\partial z}\left[\frac{1}{1-|\omega_1|^2}\left(\frac{h''_1}{h'_1}\omega'_1-\omega''_1\right)\right]\\
&-\frac32\overline{\omega_1}^{\,2}\frac{\partial }{\partial z}\left[\left(\frac{\omega'_1}{1-|\omega_1|^2}\right)^2\right]\,.
\end{align*}
Thus  we obtain
\[
\frac{\partial S_H(f_1)}{\partial z}(0)=(S(h_1))'(0)\,.
\]
\par
Similarly, we also have that for any integer $n\geq 2$,
\begin{align*}
\frac{\partial^n S_H(f_1)}{\partial z^n}&=(S(h_1))^{(n)}+\overline\omega_1\frac{\partial^n}{\partial z^n}\left[\frac{1}{1-|\omega_1|^2}\left(\frac{h''_1}{h'_1}\omega_f'-\omega''_1\right)\right]\\
&-\frac32\overline{\omega_1}^{\,2}\frac{\partial^n }{\partial z^n}\left[\left(\frac{\omega'_1}{1-|\omega_1|^2}\right)^2\right]\,.
\end{align*}
Therefore,
\[
\frac{\partial^n S_f}{\partial z^n}(0)=(Sh)^{(n)}(0)
\]
for all $n\geq 1$.
\par
By repeating the same procedure (considering now the function $f_2$ instead of $f_1$), we obtain
\[
\frac{\partial^n S_H(f_2)}{\partial z^n}(0)=(S(h_2))^{(n)}(0)\quad \text{for all}\quad  n\geq 1\,.
\]
\par
Since we are assuming that $S_H(f_1)=S_H(f_2)$, we get the identities
\[
(S(h_1))^{(n)}(0)=(S(h_2))^{(n)}(0)
\]
 for all positive integer $n$. This proves i).
\par\smallskip
To prove \eqref{eq-prop-aux1}, we are going to show that the Taylor coefficients of the (analytic) functions
\[
\omega'_1(0)\left(\frac{h''_1}{h'_1}\omega'_1-\omega''_1\right)\quad\text{and}\quad \omega'_2(0)\left(\frac{h''_2}{h'_2}\omega'_2-\omega''_2\right)
\]
coincide.
\par
Using that $S(h_1)=S(h_2)$ and \eqref{eq-prop-equalSh}, we have
\begin{align}\label{eq-prop-equalSh2}
\nonumber \frac{\overline \omega_1}{1-|\omega_1|^2}& \left(\frac{h''_1}{h'_1}\,\omega'_1-\omega''_1\right)
-\frac 32\left(\frac{\omega'_1\,\overline
\omega_1}{1-|\omega_1|^2}\right)^2\\
&  =\frac{\overline \omega_2}{1-|\omega_2|^2} \left(\frac{h''_2}{h'_2}\,\omega'_2-\omega''_2\right)
-\frac 32\left(\frac{\omega'_2\,\overline\omega_2}{1-|\omega_2|^2}\right)^2\,.
\end{align}
Taking derivatives with respect to $\overline z$ in both sides of the previous equation, we get
\begin{align}
& \nonumber \frac{\overline{\omega'_1}}{(1-|\omega_1|^2)^2} \left(\frac{h''_1}{h'_1}\,\omega'_1-\omega''_1\right)
-\overline{\omega_1} \left(\, \frac{3\omega'_1|\omega'_1|^2}{(1-|\omega_1|^2)^3}\right)\\
& \label{eq-Lemma2 1}  =\frac{\overline{\omega'_2}}{(1-|\omega_2|^2)^2} \left(\frac{h''_2}{h'_2}\,\omega'_2-\omega''_2\right)
-\overline{\omega_2} \left(\, \frac{3\omega'_2|\omega'_2|^2}{(1-|\omega_2|^2)^3}\right)\,.
\end{align}
This implies, evaluating at $z=0$ and using also \eqref{eq-normalizationw} and \eqref{eq-normalizationw'},
\[
\omega'_1(0)\left(\frac{h''_1}{h'_1}\omega'_1-\omega''_1\right)(0)= \omega'_2(0)\left(\frac{h''_2}{h'_2}\omega'_2-\omega''_2\right)(0)\,.
\]
Now note that
\begin{equation*}
\frac{\partial}{\partial z} \left(\frac{\overline{\omega'}}{(1-|\omega|^2)^2}\right)= \overline \omega \left(\frac{2|\omega'|^2}{(1-|\omega|^2)^3}\right)\,.
\end{equation*}
Hence, taking the derivatives with respect to $z$ of the functions in \eqref{eq-Lemma2 1} gives
\begin{align}
\nonumber \frac{\overline{\omega'_1}}{(1-|\omega_1|^2)^2} &  \left(\frac{h''_1}{h'_1}\,\omega'_1-\omega''_1\right)'
-\overline{\omega_1} \varphi_1\\ &=
\label{eq-lemma 2 3}\frac{\overline{\omega'_2}}{(1-|\omega_2|^2)^2} \left(\frac{h''_2}{h'_2}\,\omega'_2-\omega''_2\right)'
-\overline{\omega_2} \varphi_2
\end{align}
for appropriate (smooth functions) $\varphi_1$ and $\varphi_2$. This gives (since $\omega_1(0)=\omega_2(0)=0$)
\[
\omega'_1(0)\left(\frac{h''_1}{h'_1}\omega'_1-\omega''_1\right)'(0)= \omega'_2(0)\left(\frac{h''_2}{h'_2}\omega'_2-\omega''_2\right)'(0)\,.
\]
By taking successive derivatives with respect to $z$ in \eqref{eq-lemma 2 3} and evaluating at the origin, we get the desired result.
\par
To check \eqref{eq-prop-aux2}, let us introduce the notation
\[
\Phi_1=\frac{h''_1}{h'_1}\omega'_1-\omega''_1\quad\text{and}\quad \Phi_2=\frac{h''_2}{h'_2}\omega'_2-\omega''_2\,,
\]
so that, after taking complex conjugates, \eqref{eq-prop-equalSh2} becomes
\[
\frac{\overline{\Phi_1}}{1-|\omega_1|^2}\omega_1-\frac 32 \left(\frac{\overline{\omega'_1}}{1-|\omega_1|^2}\right)^2(\omega_1)^2
=\frac{\overline{\Phi_2}}{1-|\omega_2|^2}\omega_2-\frac 32 \left(\frac{\overline{\omega'_2}}{1-|\omega_2|^2}\right)^2(\omega_2)^2\,.
\]
Taking derivatives with respect to $z$ of the functions in the previous identity and using that
\[
\frac{\partial}{\partial z} \frac{1}{1-|\omega|^2}=\overline{\omega}\frac{\omega'}{(1-|\omega|^2)^2}\,,
\]
we obtain
\begin{align*}
\frac{\overline{\Phi_1}}{1-|\omega_1|^2}\omega'_1
&-\frac 32 \left(\frac{\overline{\omega'_1}}{1-|\omega_1|^2}\right)^2\left((\omega_1)^2\right)'+\overline{\omega_1} \psi_1 \\ &=\frac{\overline{\Phi_2}}{1-|\omega_2|^2}\omega'_2
-\frac 32 \left(\frac{\overline{\omega'_2}}{1-|\omega_2|^2}\right)^2\left((\omega_2)^2\right)'+\overline{\omega_2} \psi_2
\end{align*}
for appropriate smooth functions $\psi_1$ and $\psi_2$. Therefore, in view of the previous identity, we can argue as above (calculating successive derivatives with respect to $z$ and evaluating at the origin) to show that the analytic functions in \eqref{eq-prop-aux2} have equal Taylor coefficients, so that they are equal. We omit the details.
\end{pf}

Some remarks are now in order. Notice that the condition that $S(h_1)=S(h_2)$, in addition to the normalizations \eqref{eq-normalizationh}, show that there exists a constant $a_0$, say, such that
\begin{equation}\label{eq-afterprop h}
h_1(z)=\frac{h_2(z)}{1+a_0h_2(z)}\,,\quad z\in\mathbb{D}\,.
\end{equation}
A straightforward calculation shows that $a_0=(h''_2(0)-h''_1(0))/2$.
\par\smallskip
We have not been able to use all the information contained in the equations \eqref{eq-prop-aux1}, \eqref{eq-prop-aux2}, and \eqref{eq-afterprop h} to get the solution to the problem we are to solve directly. However, under one single extra assumption on the values of the derivatives of the dilatations $\omega_1$ and $\omega_2$ at the origin, it is possible to identify the relationship between the harmonic mappings $f_1$ and $f_2$ (normalized as in the previous proposition) such that $S_H(f_1)=S_H(f_2)$.

\begin{prop}\label{prop-equalfucn}
Assume that the harmonic mappings $f_1=\overline{g_1}+h_1$ and $f_2=\overline{g_2}+h_2$ satisfy the hypotheses in Proposition~\ref{prop-schw h}. If, in addition, $\omega'_1(0)=\omega'_2(0)$, then $f_1=f_2$.
\end{prop}

\begin{pf}
The hypotheses show that, by \eqref{eq-prop-aux1},
\begin{equation*}
\frac{h''_1}{h'_1}\, \omega'_1-\omega''_1=\frac{h''_2}{h'_2}\, \omega'_2-\omega''_2
\end{equation*}
and hence
\begin{equation*}
\frac{h''_1(0)}{h'_1(0)}\, \omega'_1(0)-\omega''_1(0)=\frac{h''_2(0)}{h'_2(0)}\, \omega'_2(0)-\omega''_2(0)\,.
\end{equation*}
Therefore, from \eqref{eq-prop-aux2} we get
\[
\overline{\left(\frac{h''_1(0)}{h'_1(0)}\, \omega'_1(0)-\omega''_1(0)\right)} \left(\omega_1-\omega_2\right)=\frac 32 (\omega'_1(0))^2 \left((\omega_1)^2-(\omega_2)^2\right)
\]
and we conclude that either $\omega_1=\omega_2$ or
\[
\omega_1+\omega_2=\frac{2}{3(\omega'_1(0))^2}\overline{\left(\frac{h''_1(0)}{h'_1(0)}\, \omega'_1(0)-\omega''_1(0)\right)} \,.
\]
The latter case gives
\[
\omega'_1(0)=-\omega'_2(0)\,,
\]
which is not possible since we are assuming that the normalization \eqref{eq-normalizationw'} holds, so that $\omega'_1(0)$ and $\omega'_2(0)$ are both positive real numbers. Hence, $\omega_1\equiv\omega_2$ and, by \eqref{eq-prop-aux1}, we conclude $h''_1(0)=h''_2(0)$. This shows, by \eqref{eq-afterprop h}, that $h_1=h_2$. Finally, we have also
\[
g'_1=\omega_1 h'_1=\omega_2 h'_2=g'_2\,,
\]
$g_1(0)=g_2(0)$. This gives $g_1=g_2$ and we get the desired identity $f_1=f_2$.
\end{pf}
\par\smallskip
From now on, all our efforts will be devoted to show that if $f_1$ and $f_2$ are as in Proposition~\ref{prop-schw h}, the identity $\omega'_1(0)=\omega'_2(0)$ holds.

\subsection{An equivalent formulation of the problem}\label{ssec-equivalentformulation}

Given an orien\-ta\-tion-preserving harmonic mapping $f=\overline{g}+h$ in the unit disk with dilatation $\omega$ and normalized as in \eqref{eq-normalizationh}, \eqref{eq-normalizationw}, and \eqref{eq-normalizationw'}, and $w\in\mathbb{D}$, we can consider the transformation defined by \eqref{eq-aff1} to get a new function
\[
\widehat{F_w}=\overline{\widehat{G_w}}+\widehat{H_w}=\overline{\left(\frac{g\circ\varphi_w-g(w)}{\overline{h'(w)}(1-|w|^2)}\right)}+\frac{h\circ\varphi_w-h(w)}{h'(w)(1-|w|^2)}\,,
\]
where $\varphi_w$ is the automorphism of the unit disk given by \eqref{eq-automorphism}.
\par
The anti-analytic rotation $R_\mu(\widehat{F_w})=F_w$ (defined as in \eqref{eq-antianalyticrotation}) with $\mu=h'(w)/\overline{h'(w)}$  gives
\[
F_w=H_w+\overline{G_w}=\frac{h\circ\varphi_w-h(w)}{h'(w)(1-|w|^2)}+\overline{\left(\frac{g\circ\varphi_w-g(w)}{h'(w)(1-|w|^2)}\right)}\,.
\]
This new function $F_w$ satisfies the normalization \eqref{eq-normalizationh}. However, the dilatation $\alpha_w$ of $F_w$ equals $\alpha_w=\omega_f\circ\varphi_w$, which clearly does not fix the origin necessarily, so that \eqref{eq-normalizationw} might not be satisfied. In this case we can consider the affine mapping $A(z)=(z-\overline{\alpha_w(0)\, z})/(1-|\alpha_w(0)|^2)$  and the composition $A\circ F_w$ to get the function $\widehat{f_w}=\overline{\widehat{g_w}}+\widehat{h_w}$ with dilatation $\widehat{\omega_w}$, where
\[
\widehat{h_w}=\frac{1}{1-|\omega(w)|^2} \left(\frac{h\circ\varphi_w-h(w)}{h'(w)(1-|w|^2)}-\overline{\omega(w)} \left(\frac{g\circ\varphi_w-g(w)}{h'(w)(1-|w|^2)}\right) \right)
\]
and $\widehat{\omega_w}=\varphi_{-\omega(w)}\circ\omega_f\circ\varphi_w$.
\par
Note that $\widehat{\omega_w}(0)=0$, so that $\widehat{f_w}$ satisfies both \eqref{eq-normalizationh} and \eqref{eq-normalizationw}.  The further transformation $f_w=R_{\mu_w}(\widehat{f_w})$, where $\mu_w=\overline{\widehat{\omega'_w(0)}}/|\widehat{\omega'_w(0)}|$ produces a function of the form
\[
f_w=\overline{g_w}+h_w\,,
\]
with $h_w=\widehat{h_w}$, normalized as in \eqref{eq-normalizationh}, \eqref{eq-normalizationw}, and \eqref{eq-normalizationw'}, and (by the chain rule and the invariance of the operator $S_H$ under affine transformations and anti-analytic rotations) such that
\begin{equation}\label{eq-identitySh}
S_H(f_w)(z)=S_H(f)(\varphi_w(z))(\varphi'_w(z))^2\,,\quad z\in\mathbb{D}\,.
\end{equation}
Note that for all such $z$,
\begin{align}\label{eq-schwarzianh}
\nonumber S(h_w)(z)&=S\left(h\circ\varphi_w-\overline{\omega(w)}(g\circ\varphi_w)\right)(z)
\\&=S(h-\overline{\omega(w)}g)(\varphi_w(z))\cdot (\varphi'_w(z))^2\,.
\end{align}
\par\medskip
We can now state our next theorem.
\begin{theorem}\label{thm-mainaux}
Let $f_1=\overline{g_1}+h_1$ and $f_2=\overline{g_2}+h_2$ be orientation-preserving harmonic mappings in the unit disk with non-constant dilatations $\omega_1=g'_1/h'_1$ and $\omega_2=g'_2/h'_2$, respectively. Assume that both $f_1$ and $f_2$ are normalized as in \eqref{eq-normalizationh}, \eqref{eq-normalizationw}, and \eqref{eq-normalizationw'}. Then $S_H(f_1)=S_H(f_2)$ if and only if for all $w$ in the unit disk
\begin{equation}\label{eq-thm-schwarzianh1}
S(h_1-\overline{\omega_1(w)}g_1)=S(h_2-\overline{\omega_2(w)}g_2)\,.
\end{equation}
\end{theorem}
\begin{pf}
Suppose first that \eqref{eq-thm-schwarzianh1} holds. Then, in particular, we have that for all $w\in\mathbb{D}$
\[
S(h_1-\overline{\omega_1(w)}g_1)(w)=S(h_2-\overline{\omega_2(w)}g_2)(w)\,,
\]
which readily gives, by \eqref{eq-Schwarzian2}, $S_H(f_1)=S_H(f_2)$.
\par\smallskip
To prove the necessity of the assertion in the theorem, we argue as before to produce, for each $w$ in the unit disk the functions $(f_1)_w=\overline{(g_1)_w}+(h_1)_w$ and $(f_2)_w=\overline{(g_2)_w}+(h_2)_w$, that will satisfy \eqref{eq-normalizationh}, \eqref{eq-normalizationw}, and \eqref{eq-normalizationw'}. Moreover, due to the identity \eqref{eq-identitySh} applied, respectively, to $(f_1)_w$ and $(f_2)_w$, and since $S_H(f_1)=S_H(f_2)$ we obtain, for all $z\in\mathbb{D}$,
\begin{align*}
S_H((f_1)_w)(z)& =S_H(f_1)(\varphi_w(z)) (\varphi'_w(z))^2\\
& =S_H(f_2)(\varphi_w(z)) (\varphi'_w(z))^2=S_H((f_2)_w)(z)
\end{align*}
and hence, by Proposition~\ref{prop-schw h} i), we have $S_H((h_1)_w)=S_H((h_2)_w)$. In other words, according to \eqref{eq-schwarzianh}, we have
\begin{equation*}
S(h_1-\overline{\omega_1(w)}g_1)(\varphi_w(z))\cdot (\varphi'_w(z))^2=S(h_1-\overline{\omega_1(w)}g_1)(\varphi_w(z))\cdot (\varphi'_w(z))^2\,,
\end{equation*}
which is equivalent to \eqref{eq-thm-schwarzianh1} since $\varphi_w$ is an automorphism of $\mathbb{D}$.
\end{pf}
\par\smallskip
The following corollary will provide us with another important tool to get the solution to the problem considered in this section.
\begin{cor}
Let $f_1$ and $f_2$ satisfy the hypotheses in the previous theorem. If we assume that  $S_H(f_1) = S_H(f_2)$, then the identity
\begin{equation}\label{eq-thm-schwarzianh2}
h_1(z)-\overline{\omega_1(z)}g_1(z)=\frac{h_2(z)-\overline{\omega_2(z)}g_2(z)}{1+\left(a_0+\delta(z)\right)\left(h_2(z)-\overline{\omega_2(z)}g_2(z)\right)}\,,
\end{equation}
where
\begin{equation}\label{eq-delta}
a_0=\frac{h''_2(0)-h''_1(0)}{2}\quad\text{and}\quad \delta(z)=\frac{\omega'_1(0)}{2}\overline{\omega_1(z)}-\frac{\omega'_2(0)}{2}\overline{\omega_2(z)}\,,
\end{equation}
holds for all  $z$ in $\mathbb{D}$.
\end{cor}
\begin{pf}
According to the previous theorem, for a complex number $w\in\mathbb{D}$, the analytic functions $\gamma_1=h_1-\overline{\omega_1(w)}g_1$ and $\gamma_2=h_2-\overline{\omega_2(w)}g_2$ have equal (analytic) Schwarzian derivative. Therefore, as mentioned in the introduction, there exists a M\"{o}bius transformation $T_w$ as in \eqref{eq-Mobius} (where the coefficients can depend on $w$) such that $\gamma_1=T_w\circ\gamma_2$. But since $\gamma_1(0)=\gamma_2(0)=0$ and $\gamma'_1(0)=\gamma'_2(0)=1$, we have that, necessarily,
\[
T_w(z)=\frac{z}{1+a_w z}\,.
\]
That is,
\[
h_1(z)-\overline{\omega_1(w)}g_1(z)=\frac{h_2(z)-\overline{\omega_2(w)}g_2(z)}{1+a_w\left(h_2(z)-\overline{\omega_2(w)}g_2(z)\right)}\,.
\]
Taking two successive derivatives with respect to $z$ in both sides of the previous equation, evaluating at $z=0$, and bearing in mind that \eqref{eq-normalizationw} and \eqref{eq-normalizationw'} are satisfied by both $\omega_1$ and $\omega_2$, we get
\[
a_w= \frac{h''_2(0)-h''_1(0)}{2}+\frac{\omega'_1(0)}{2}\overline{\omega_1(w)}-\frac{\omega'_2(0)}{2}\overline{\omega_2(w)}\,.
\]
Setting $w=z$, we obtain the desired result due.
\end{pf}
\subsection{Important lemmas}\label{ssec-prel} In order to simplify the exposition, let us agree with the following notation already used in the proof of Proposition~\ref{prop-schw h}. We use $\Phi_1$ and $\Phi_2$ to denote the analytic functions in the unit disk defined by
\begin{equation}\label{eq-notation}
\Phi_1= \frac{h''_1}{h'_1}\, \omega'_1-\omega''_1\quad\text{and}\quad \Phi_2= \frac{h''_2}{h'_2}\, \omega'_2-\omega''_2\,,
\end{equation}
so that \eqref{eq-prop-aux1} and \eqref{eq-prop-aux2} become
\begin{equation}\label{eq-afterprop-aux1}
\omega'_1(0) \Phi_1=\omega'_2(0) \Phi_2
\end{equation}
and
\begin{equation}\label{eq-afterprop-aux2}
\overline{\Phi_1(0)} \omega_1-\frac 32 (\omega'_1(0))^2\, (\omega_1)^2=\overline{\Phi_2(0)} \omega_2-\frac 32 (\omega'_2(0))^2\, (\omega_2)^2\,,
\end{equation}
respectively.
\par\smallskip
\begin{lem}\label{lem-Phi}
Let $f_1=\overline{g_1}+h_1$ and $f_2=\overline{g_2}+h_2$ be orientation-preserving harmonic mappings in the unit disk with non-constant dilatations $\omega_1$ and $\omega_2$, respectively. Assume that both $f_1$ and $f_2$ are normalized as in \eqref{eq-normalizationh}, \eqref{eq-normalizationw}, and \eqref{eq-normalizationw'}. Suppose that $S_H(f_1)=S_H(f_2)$. If either $\Phi_1(0)=0$ or $\Phi_2(0)=0$, where $\Phi_1$ and $\Phi_2$ are as in \eqref{eq-notation}, then $\omega'_1(0)=\omega'_2(0)$.
\end{lem}
\begin{pf}
If  $\Phi_1(0)=0$, by \eqref{eq-afterprop-aux1} and the normalization \eqref{eq-normalizationw'}, $\Phi_2(0)=0$ as well, so that \eqref{eq-afterprop-aux2} becomes
\[
(\omega'_1(0))^2\, (\omega_1)^2=(\omega'_2(0))^2\, (\omega_2)^2\,.
\]
By taking two successive derivatives in both sides of the previous equation and evaluating at the origin, we get the identity $\omega'_1(0)=\omega'_2(0)$. It is obvious that the same result holds if we suppose $\Phi_2(0)=0$.
\end{pf}
\par\smallskip
The next result will be also used later.
\begin{lem}\label{lem-dilatations}
Let $\omega_1$ and $\omega_2$ be the dilatations of two harmonic functions $f_1$ and $f_2$ as in the previous lemma. Assume that some of the following relations hold in some open set $\Delta$ contained in the unit disk for certain constants $k, l,$ and $m$ in $\mathbb{C}$.
\begin{itemize}
\item[i)] $\omega_2=k\omega_1\,.$
\item[ii)] $\omega_1\omega_2=k(\omega_2)^2+l\omega_2+m\omega_1\,.$
\end{itemize}
Then $\omega'_1(0)=\omega'_2(0)$.
\par\smallskip
Moreover, such functions $\omega_1$ and $\omega_2$ cannot satisfy a relation of the form
\begin{equation}\label{eq-lemdilatfinal}
n(\omega_2)^2=k\omega_2+l\omega_1+m\,,
\end{equation}
unless $n=0$.
\end{lem}
\begin{pf}
Notice that by the uniqueness principle for analytic functions, if any of the relations considered in this lemma hold in $\Delta$, they indeed hold in the whole unit disk.
\par
Assume first that i) holds. Then, necessarily, $k\neq 0$ since, otherwise \eqref{eq-normalizationw'} would not be satisfied by $\omega_2$. In fact, by taking derivatives in both sides of the equation and evaluating at the origin we obtain $k=\omega'_2(0)/\omega'_1(0)$. Hence, the substitution $\omega_2=k\omega_1$ in \eqref{eq-afterprop-aux2} and \eqref{eq-afterprop-aux1} give
\begin{equation*}
(\omega'_1(0))^2 (\omega_1)^2 =\frac{(\omega'_2(0))^4}{(\omega'_1(0))^2} (\omega_1)^2\,,
\end{equation*}
which implies (taking two successive derivatives and evaluating at the origin) $\omega'_1(0)=\omega'_2(0)$.
\par
Let us suppose now that for some constants $k, l,$ and $m$,
\begin{equation}\label{eq-ii.}
\omega_1\omega_2=k(\omega_2)^2+l\omega_2+m\omega_1\,.
\end{equation}
This gives, in particular (after taking derivatives with respect to $z$ and using the normalizations \eqref{eq-normalizationw}), that $0=l\omega'_2(0)+m\omega'_1(0)$. Moreover, if $l=m=0$, we get $\omega_1=k\omega_2$, which implies $\omega'_1(0)=\omega'_2(0)$ by i) and, also, if $(l,m)\neq (0,0)$, we can take derivatives in \eqref{eq-ii.} and evaluate at the origin to get that either $\omega'_1(0)$ or $\omega'_2(0)$ are zero, which is in contradiction with \eqref{eq-normalizationw'}. So that we can assume that both $l$ and $m$ are different from zero and, in fact, we obtain
\begin{equation}\label{eq-L}
l=-m\omega'_1(0)/\omega'_2(0)\,.
\end{equation}
Also, we can re-write \eqref{eq-ii.} as
\[
\omega_1=\frac{k(\omega_2)^2+l\omega_2}{\omega_2-m}
\]
and substitute this expression in \eqref{eq-afterprop-aux2} to get (after multiplying by $(\omega_2-m)^2$, re-arranging the resulting terms, using the fact that $l=-m\omega'_1(0)/\omega'_2(0)$, and \eqref{eq-afterprop-aux1})
\begin{align*}
\frac 32&\left(k^2(\omega'_1(0))^2-(\omega'_2(0))^2\right)(\omega_2)^2\\
&+\left(\overline{\Phi_2(0)}+3m(\omega'_2(0))^2-k\overline{\Phi_1(0)}+3kl(\omega'_1(0))^2\right)\omega_2\\
& +\left(\frac 32 l^2(\omega'_1(0))^2+(mk-l)\overline{\Phi_1(0)}-2m\overline{\Phi_2(0)}-\frac 32 m^2(\omega'_2(0))^2\right)\equiv 0\,.
\end{align*}
Hence, we have that the coefficients in the previous equation are identically zero. That is,
\begin{equation}\label{eq-K}
k^2=\left(\frac{\omega'_2(0)}{\omega'_1(0)}\right)^2\,,
\end{equation}
\begin{align}\label{eq-aux1}
\nonumber k\overline{\Phi_1(0)}-\overline{\Phi_2(0)}&=3m(\omega'_2(0))^2+3kl(\omega'_1(0))^2\\
&= 3m(\omega'_2(0))^2-3km\frac{(\omega'_1(0))^3}{\omega'_2(0)}
\end{align}
(by \eqref{eq-L}), and, using again \eqref{eq-L} and also \eqref{eq-afterprop-aux1},
\begin{equation}\label{eq-aux22}
k\overline{\Phi_1(0)}-\overline{\Phi_2(0)}=\frac 32 m \frac{(\omega'_2(0))^4-(\omega'_1(0))^4}{(\omega'_2(0))^2}\,.
\end{equation}
A combination of \eqref{eq-aux1} and \eqref{eq-aux22} gives the identity
\[
4k^2(\omega'_1(0))^6 (\omega'_2(0))^2=(\omega'_1(0))^8+(\omega'_2(0))^8+2(\omega'_1(0))^4(\omega'_2(0))^4\,,
\]
so that, in view of \eqref{eq-K}, we finally get $\omega'_1(0)=\omega'_2(0)$.
\par\smallskip
Finally, to prove the last assertion of the theorem, suppose that \eqref{eq-lemdilatfinal} is satisfied for some $n\neq 0$ so that, we can assume, without loss of generality that $n=1$ to get a relation of the form
\[
(\omega_2)^2=k\omega_2+l\omega_1+m\,,\quad k, l, m\in\mathbb{C}\,.
\]
Then, on the one hand, $m=0$ since $\omega_1(0)=\omega_2(0)=0$. Also, $l\neq 0$ since $\omega_2$ is not identically constant. On the other hand, $k\neq 0$ since $\omega'_1(0)\neq 0$. Hence, we can re-write the previous equation as
\[
\omega_1=\frac{1}{l}(\omega_2)^2-\frac{k}{l}\ \omega_2\,,\quad k, l\in\mathbb{C}\setminus\{0\}\,.
\]
But then, after using this identity in \eqref{eq-afterprop-aux1} and rearranging the resulting terms, we get
\[
a_4(\omega_2)^4+a_3(\omega_2)^3+a_2(\omega_2)^2+a_1\omega_2\equiv 0
\]
(for certain coefficients $a_1, a_2, a_3$, and $a_4$). Therefore, since $\omega_2(\mathbb{D})$ is an open set we conclude all the coefficients in the previous equation must be equal to zero. In particular,
\[
a_3=-\frac{3(\omega'_1(0))^2k}{l^2}=0\,,
\]
which implies, since $\omega'_1(0)\neq 0$,  that $k=0$. This is a contradiction.
\end{pf}
\par\smallskip
The last result in this section is related to the following analytic functions that will be important in the proof of our main theorem in the next Section~\ref{ssec-equaldilat}. They are defined in terms of the functions in the canonical decomposition of $f_1=\overline{g_1}+h_1$ and $f_2=\overline{g_2}+h_2$, the value of the derivative of their dilatations $\omega_1$ and $\omega_2$, respectively, at the origin, and the constant $a_0$ in \eqref{eq-delta}. Concretely, define
\begin{equation}\label{deq-defn1}
\varphi_1=g_2-a_0h_1g_2-\frac{\omega'_2(0)}{2}h_1h_2\,,\quad \varphi_2=\frac{\omega'_2(0)}{2}h_1g_2\,,
\end{equation}
\begin{equation}\label{deq-defn2}
\varphi_3=\frac{\omega'_1(0)}{2}g_1g_2\,,\quad \varphi_4=-\frac{\omega'_2(0)}{2}g_1g_2\,,
\end{equation}
\begin{equation}\label{deq-defn3}
\varphi_5=a_0g_1g_2-\frac{\omega'_1(0)}{2}h_1g_2+\frac{\omega'_2(0)}{2}g_1h_2\,,\quad \varphi_6=-\frac{\omega'_1(0)}{2}g_1h_2\,,
\end{equation}
\begin{equation}\label{deq-defn4}
\text{and}\quad  \varphi_7=\frac{\omega'_1(0)}{2}h_1h_2-a_0g_1h_2-g_1\,.
\end{equation}
\par
Notice that, since we are assuming that the functions $f_1$ and $f_2$ are orientation-preserving (so that both $h_1$ and $h_2$ are analytic and locally univalent) and that the dilatations $\omega_1=g'_1/h'_1$ and $\omega_2=g'_2/h'_2$ are not constant and satisfy \eqref{eq-normalizationw'}, we have that neither $\varphi_2$, $\varphi_3$, $\varphi_4$, nor $\varphi_6$ can be identically zero. At this point, and due to the remarks we made before, we can also prove that $\varphi_1$, $\varphi_5$, and $\varphi_7$ cannot be identically zero either.
\begin{lem}\label{lem-varphi}
Under the hypotheses considered in Lemma~\ref{lem-Phi} and the additional condition $\omega'_1(0)\neq\omega'_2(0)$, we have that neither $\varphi_1$, $\varphi_5$, nor $\varphi_7$ are identically zero.
\end{lem}
\begin{pf}
Assume that $\varphi_1\equiv 0$. Then we have, in view of \eqref{eq-afterprop h},
\[
g_2=\frac{\omega'_2(0)}{2}\frac{h_1}{1-a_0h_1}h_2=\frac{\omega'_2(0)}{2}(h_2)^2\,.
\]
Hence, taking successive derivatives and using \eqref{eq-normalizationh}, we get that $g'''_2(0)=3\omega'_2(0)h''_2(0)$. On the other hand, $g'_2=\omega_2 h'_2$, so that taking derivatives again, we obtain the equation
\begin{equation}\label{eq-g'''}
g'''_2(0)=\omega''_2(0)+2\omega'_2(0)h''_2(0)\,.
\end{equation}
That is,
\[
3\omega'_2(0)h''_2(0)=\omega''_2(0)+2\omega'_2(0)h''_2(0)\,,
\]
or equivalently, $\Phi_2(0)=0$, where $\Phi_2$ is the function defined in \eqref{eq-notation}. But then, by Lemma~\ref{lem-Phi}, we obtain that $\omega'_1(0)=\omega'_2(0)$, which is in contradiction with our hypotheses.
\par
The proof that $\varphi_7$ cannot be identically zero either is completely analogous. We omit the details.
\par\smallskip
Let us suppose now that $\varphi_5\equiv 0$ and recall that, as in \eqref{eq-delta}, $a_0=(h''_2(0)-h''_1(0))/2$. Then
\[
\lim_{z\to 0} a_0 \frac{g_1g_2}{z^4}=\lim_{z\to 0} \frac{\frac{\omega'_1(0)}{2}h_1g_2-\frac{\omega'_2(0)}{2}h_2g_1}{z^4}\,.
\]
A straightforward calculation, which uses that $g'_i=\omega_ih'_i$, $i=1,2$, and that \eqref{eq-normalizationh} and \eqref{eq-normalizationw} are satisfied for all the functions involved, shows that
\[
\lim_{z\to 0} a_0 \frac{g_1g_2}{z^4}=a_0\frac{\omega'_1(0)\omega'_2(0)}{4}\,.
\]
With some for effort and using \eqref{eq-g'''} and the analogous identity $g'''_1(0)=\omega''_1(0)+2\omega'_1(0)h''_1(0)$, we obtain
\[
\lim_{z\to 0} \frac{\frac{\omega'_1(0)}{2}h_1g_2-\frac{\omega'_2(0)}{2}h_2g_1}{z^4}=\frac{\omega'_1(0)\omega''_2(0)-\omega''_1(0)\omega'_2(0)+a_0\omega'_1(0)\omega'_2(0)}{12}\,.
\]
Hence, we conclude
\[
\omega'_1(0)\left(\omega'_2(0)h''_2(0)-\omega''_2(0)\right)=\omega'_2(0)\left(\omega'_1(0)h''_1(0)-\omega''_1(0)\right)
\]
or equivalently, using also \eqref{eq-afterprop-aux1},
\[
\omega'_1(0) \Phi_2(0)=\omega'_2(0) \Phi_1(0)= \omega'_2(0) \frac{\omega'_2(0)}{\omega'_1(0)}\Phi_2(0)\,.
\]
This gives (using Lemma~\ref{lem-Phi} if needed) that $\omega'_1(0)=\omega'_2(0)$. This contradiction ends the proof of this lemma.
\end{pf}
\subsection{Main theorem} \label{ssec-equaldilat}
Now we have all the tools to prove the main result in this section. We are aware of the fact that the proof of this theorem is very technical: there are numerous different equations that must be combined appropriately to get the desired conclusion. With the hope to make our arguments understandable and this paper self-contained, we have decided to include enough details in the different steps in our approach to the proof.
\begin{theorem}\label{thm-main}
Let $f_1=\overline{g_1}+h_1$ and $f_2=\overline{g_2}+h_2$ be orientation-preserving harmonic mappings in the unit disk with non-constant dilatations $\omega_1$ and $\omega_2$, respectively. Assume that both $f_1$ and $f_2$ are normalized as in \eqref{eq-normalizationh}, \eqref{eq-normalizationw}, and \eqref{eq-normalizationw'}. Suppose that $S_H(f_1)=S_H(f_2)$. Then $\omega'_1(0)=\omega'_2(0)$.
\end{theorem}
\begin{pf}
Let us assume, in order to get a contradiction, that $\omega'_1(0)\neq\omega'_2(0)$.
\smallskip
\par
After multiplying both sides of \eqref{eq-thm-schwarzianh2} by $(1+(a_0+\delta)(h_2-\overline{\omega_2}g_2))$, bearing in mind the definition of the function $\delta$ in \eqref{eq-delta}, using that from \eqref{eq-afterprop h}, we have $h_2-h_1-a_0h_1h_2=0$, and rearranging the resulting terms, we obtain the relation
\begin{equation}\label{eq-1}
0=\varphi_1 \overline{B_1}+\varphi_2 \overline{B_2}+\varphi_3 \overline{B_3}+\varphi_4 \overline{B_4}+\varphi_5 \overline{B_5}+\varphi_6 \overline{B_6}+\varphi_7 \overline{B_7}\,,
\end{equation}
where for $i=1,\ldots, 7$, the functions $\varphi_i$ are defined by \eqref{deq-defn1}, \eqref{deq-defn2}, \eqref{deq-defn3}, and \eqref{deq-defn4}, and the analytic functions $B_i$ in the unit disk are
\begin{align}\label{deq-defn5}
\nonumber B_1&=\omega_2\,,\quad B_2=(\omega_2)^2\,, \quad B_3=(\omega_1)^2\omega_2\,,\quad B_4=\omega_1(\omega_2)^2\,,\\
& \quad B_5=\omega_1\omega_2\,,\quad B_6=(\omega_1)^2\,,\quad\text{and}\quad   B_7=\omega_1\,.
\end{align}

\par\smallskip
Now, notice that the fact that $\omega'_1(0)>0$ (which implies that $\omega_1$ is univalent in a disk centered at the origin) shows that for all $i=1,\ldots, 6$, the quotients $B_i/B_7$ are analytic functions on a neighborhood of the origin.  Hence, if we divide by $\overline{B_7}$ in \eqref{eq-1} and make the derivatives with respect to $\overline z$ in both sides of the resulting equation, we obtain
\begin{equation}\label{eq-2}
0=\varphi_1 \overline{C_1}+\varphi_2 \overline{C_2}+\varphi_3 \overline{C_3}+\varphi_4 \overline{C_4}+\varphi_5 \overline{C_5}+\varphi_6 \overline{C_6}\,,
\end{equation}
where $C_i=(B_i/B_7)'$, $i=1,\ldots, 6$. Notice that $C_6=\omega'_1$, which is supposed to be different from zero near $z=0$. Hence, we can again divide by the complex conjugate of this function $C_6$ in the previous equation and take derivatives with respect to $\overline z$ to get
\begin{equation}\label{eq-3}
0=\varphi_1 \overline{D_1}+\varphi_2 \overline{D_2}+\varphi_3 \overline{D_3}+\varphi_4 \overline{D_4}+\varphi_5 \overline{D_5}\,,
\end{equation}
where, now, $D_i=(C_i/C_6)'$, $i=1, \ldots 5$ and, in particular, $D_5=(\omega'_2/\omega'_1)'$.

\par
Let us suppose that $D_5$ is identically zero in a neighborhood of the origin $\Delta$, say. Then, we have (after making the integration in both sides of the resulting equation $D_5=0$) that $\omega'_2=k\omega'_1$ for some constant $k\in\mathbb{C}$. This is not possible unless $\omega'_1(0)=\omega'_2(0)$, by Lemma~\ref{lem-dilatations}. Hence,  there must be a disk $\Delta_1=\Delta_1(z_0, r) \subset \Delta$ centered at $z_0\in\Delta$ and with radius $r>0$ such that the analytic function $D_5$ satisfies $D_5(z)\neq 0$ for all $z\in\Delta_1$. And hence, for all points in this disk, it makes sense do divide out by $\overline{D_5}$ in \eqref{eq-3} and take derivatives with respect to $\overline z$ to obtain
\begin{equation}\label{eq-4}
0=\varphi_1 \overline{E_1}+\varphi_2 \overline{E_2}+\varphi_3 \overline{E_3}+\varphi_4 \overline{E_4}\,,
\end{equation}
where $E_i=(D_i/D_5)'$, $i=1, \ldots 4$. We indeed have
\[
E_4=2\left(\frac{\left[\frac{\omega_2\omega'_2}{\omega'_1}\right]'}{\left[\frac{\omega'_2}{\omega'_1}\right]'}\right)'\,.
\]
Now we can argue as before to prove that $E_4$ cannot be identically zero in $\Delta_1$, as if this would be the case, we would obtain (after integrating the equation $E_4=0$ three successive times) a relation of the form
\[
(\omega_2)^2=k\omega_2+l\omega_1+m\,,\quad k, l, m\in\mathbb{C}\,,
\]
which is known to be impossible to be satisfied by Lemma~\ref{lem-dilatations}. This proves that $E_4$ is not identically zero in $\Delta_1$ and hence we have that there is another disk $\Delta_2\subset\Delta_1$ where $E_4$ has no zeros. We can then divide out by $\overline{E_4}$ in \eqref{eq-4}, take derivatives with respect to $\overline z$, and get the following equation that holds in this smaller disk $\Delta_2$:
\begin{equation}\label{eq-5}
0=\varphi_1 \overline{F_1}+\varphi_2 \overline{F_2}+\varphi_3 \overline{F_3}\,,
\end{equation}
where $F_i=(E_i/E_4)'$, for $i=1,2,$ and $3$.
\par
We still need to repeat the arguments used before to show that $F_3$ is not identically zero. To this end, it might be useful to point out that
\[
E_3=\left(\frac{\left[\frac{(\omega_1\omega_2)'}{\omega'_1}\right]'}{\left[\frac{\omega'_2}{\omega'_1}\right]'}\right)'\,.
\]
\par
So that let us assume that $F_3\equiv 0$. This gives, after four successive integrals the expression
\[
\omega_1\omega_2=k(\omega_2)^2+l\omega_2+m\omega_1\,,
\]
since $\omega_1(0)=\omega_2(0)=0$, and once more we have (after an application of Lemma~\ref{lem-dilatations}) that this would imply the contradiction $\omega'_1(0)=\omega'_2(0)$. This proves that $F_3$ is not identically zero in $\Delta_2$.
\par\smallskip
Therefore, going back to \eqref{eq-5}, we have (by Lemma~\ref{lem-varphi} and the hypotheses that the functions $f_1$ and $f_2$ considered are locally univalent and have non-constant dilatations) that none of the functions $\varphi_1$, $\varphi_2$, and $\varphi_3$ are identically zero. We have also proved that $F_3$ cannot be identically zero either. Moreover, it is now easy to obtain that $F_2$ satisfies the same property if we argue as follows.
\par
Suppose that $F_2\equiv 0$. Then, we get from \eqref{eq-5} that $\varphi_1\overline{F_1}+\varphi_3\overline{F_3}\equiv 0$. And hence, $F_1$ is not identically zero and we get the identity
\[
\frac{\varphi_1}{\varphi_3}=-\overline{\left(\frac{F_3}{F_1}\right)}\,,
\]
so that, in particular, being the functions $\varphi_1/\varphi_3$ and $F_3/F_1$ involved in the previous equation analytic in some disk $\Delta_2\subset\mathbb{D}$, we conclude that $\varphi_1=k\varphi_3$ for some constant $k\in\mathbb{C}$ different from zero. This identity is valid in $\Delta_2$. But again a direct application of the identity principle for analytic mappings gives that it must hold in the whole unit disk. A straightforward calculation shows
\[
\lim_{z\to 0}\frac{\varphi_3(z)}{z^3}=0\,,
\]
while a repeated use of L'Hopital's rule in addition to the identity $g'''_2(0)=\omega''_2(0)+2\omega'_2(0)h''_2(0)$ (obtained from $g'_2=\omega_2 h'_2$), and the fact that $a_0=(h''_2(0)-h''_1(0))/2$, gives
\[
\lim_{z\to 0}\frac{\varphi_1(z)}{z^3}=-\frac{\Phi_2(0)}{6}\left(=k\lim_{z\to 0}\frac{\varphi_3(z)}{z^3}\right)\,.
\]
Hence, we obtain $\Phi_2(0)=0$ which, as we know from Lemma~\ref{lem-Phi}, implies the identity $\omega'_1(0)=\omega'_2(0)$.  Therefore, we have that $\varphi_1$ is not a constant multiple of $\varphi_3$, so that $F_2$ is not identically zero and we can divide \eqref{eq-5} by $\overline{F_2}$ and take derivatives with respect to $\overline z$ to get
\[
\varphi_1\overline{\left(\frac{F_1}{F_2}\right)'}+\varphi_3\overline{\left(\frac{F_3}{F_2}\right)'}\equiv 0\,.
\]
This implies (using that $\varphi_1$ is not a constant multiple of $\varphi_3$) that both the functions
\[
\overline{\left(\frac{F_1}{F_2}\right)'}\quad\text{and} \quad \overline{\left(\frac{F_3}{F_2}\right)'}
\]
are identically zero, which gives the relations
\[
F_1=\widetilde{k}F_2 \quad\text{and} \quad F_3=a F_2\,,
\]
where $a\neq 0$ (since $F_3$ is not identically zero). Moreover, $\widetilde{k}$ must be different from zero as well since, otherwise, we would have $F_1\equiv 0$ and we could argue as before to get from \eqref{eq-5} that  $\varphi_2$ and $\varphi_3$ are linearly dependent, which is clearly absurd since $\omega_1$ and $\omega_2$ are not constant. That is, we can write
\begin{equation}\label{eq-EFES}
F_1=kF_3 \quad\text{and} \quad F_2=a F_3\,,\quad k, a\in\mathbb{C}\setminus\{0\}\,.
\end{equation}
Moreover, bearing in mind the definition of the functions $F_1$, $F_2$,  and $F_3$, we easily prove the following lemma.
\begin{lem}\label{lem-recursive}
The following identities hold for certain constants $a\neq 0$, $k\neq 0$, $b$, $c$, $d$, $e$ $l$, $m$, $n$, and $p$.
\begin{itemize}
\item[i)] $E_1=kE_3+lE_4$ and $E_2=aE_3+bE_4$.
\item[ii)] $D_1=kD_3+lD_4+mD_5$ and $D_2=aD_3+bD_4+cD_5$.
\item[iii)] $C_1=kC_3+lC_4+mC_5+nC_6$ and $C_2=aC_3+bC_4+cC_5+d C_6$.
\item[iv)] $B_1=kB_3+lB_4+mB_5+nB_6+pB_7$ and $B_2=aB_3+bB_4+cB_5+dB_6+eB_7$.
\end{itemize}
Moreover, from \emph{iv)},  we obtain
\begin{equation}\label{eq-I2}
\omega_2=k(\omega_1)^2\omega_2+l\omega_1(\omega_2)^2+m\omega_1\omega_2+n(\omega_1)^2+p\omega_1
\end{equation}
and
\begin{equation}\label{eq-I3}
(\omega_2)^2=a(\omega_1)^2\omega_2+b\omega_1(\omega_2)^2+c\omega_1\omega_2+d(\omega_1)^2\,,
\end{equation}
respectively, where
\begin{equation}\label{eq-P}
p=\omega'_2(0)/\omega'_1(0)\quad \text{and} \quad d=p^2-cp\,.
\end{equation}
\end{lem}
\begin{pf}
Bearing in mind that $F_1=(E_1/E_4)'$ and $F_3=(E_3/E_4)'$, we easily get from the identity $F_1=kF_3$, $k\neq 0$, obtained above, the relation
\[
\left(\frac{E_1}{E_4}\right)=k\left(\frac{E_3}{E_4}\right)+l\,.
\]
This gives the first identity in i).
\par
Now, from this first identity in i), and since $E_i=(D_i/D_5)'$, $i=1, \ldots, 4$, we obtain the first identity ii), after integrating the corresponding relation obtained in i). The same approach, recalling that we are assuming the normalizations \eqref{eq-normalizationw}, that  $D_i=(C_i/C_6)'$, $i=1, \ldots 5$, and that $C_i=(B_i/B_7)'$, $i=1,\ldots, 6$, can be used to prove that the first assertions in iii) and iv) hold.
\par\smallskip
Finally, to get \eqref{eq-I2}, it suffices to use the definition of the functions $B_i$, $i=1,\ldots, 7$ in \eqref{deq-defn5} and that \eqref{eq-normalizationw} is satisfied for both $\omega_1$ and $\omega_2$.
\par
The proofs of the second assertions in the different items in this lemma (that make use of the identity $F_2=a F_3$ obtained above), as well as the proof that \eqref{eq-I3} holds, are completely analogous to the ones we have presented. We leave the details to the reader.
\end{pf}
\par\smallskip
It is now needed to analyze further the equations \eqref{eq-I2} and \eqref{eq-I3} in the previous lemma.
\par\smallskip
Multiply \eqref{eq-I2} by $\overline{\Phi_2(0)}$, \eqref{eq-I3} by $3(\omega'_2(0))^2/2$, subtract the equations obtained, and use \eqref{eq-afterprop-aux1} and \eqref{eq-afterprop-aux2} to get the identity
\begin{align*}\label{eq-I4}
\left(\overline{\Phi_2(0)} k\right. & \left.-\frac 32 (\omega'_2(0))^2 a\right)\omega_1\omega_2+\left(\overline{\Phi_2(0)}l-\frac 32 (\omega'_2(0))^2b\right)(\omega_2)^2\\
 & +\left(\overline{\Phi_2(0)} m-\frac 32 (\omega'_2(0))^2c\right)\omega_2\\
& +\left(\overline{\Phi_2(0)} n-\frac 32 (\omega'_2(0))^2 d+\frac 32 (\omega'_1(0))^2 \right)\omega_1\equiv 0\,.
\end{align*}
A direct application of Lemma~\ref{lem-dilatations} and the fact that $\Phi_2(0)\neq 0$ by Lemma~\ref{lem-Phi} give
\begin{equation}\label{eq-identitycoefficients}
k=a\alpha\,,\quad l=b\alpha\,,\quad m=c\alpha\,,\quad \text{and} \quad n=d\alpha+\beta\,,
\end{equation}
where
\[
\alpha=\frac 32\left(\frac{(\omega'_2(0))^2}{\overline{\Phi_2(0)}}\right)\quad \text{and} \quad \beta=-\frac 32\left(\frac{(\omega'_1(0))^2}{\overline{\Phi_2(0)}}\right)\,.
\]
\par
Our next step is to use \eqref{eq-1}, \eqref{eq-2}, \eqref{eq-3}, \eqref{eq-4}, and \eqref{eq-EFES} to identify completely the constants in Lemma~\ref{lem-recursive}.
\par
To this end, first notice that by \eqref{eq-EFES} and \eqref{eq-5}, we have
\begin{equation}\label{eq-II1}
\varphi_1 \overline{a\alpha}+\varphi_2 \overline{a}+\varphi_3\equiv 0\,.
\end{equation}
\par
The substitution of the equations in assertion i)  in Lemma~\ref{lem-recursive} into the identity \eqref{eq-4} gives (using also \eqref{eq-identitycoefficients} and \eqref{eq-II1})
\begin{equation}\label{eq-II2}
\varphi_1 \overline{b\alpha}+\varphi_2 \overline{b} +\varphi_4\equiv 0\,.
\end{equation}
Hence, if we multiply \eqref{eq-II1} by $\omega'_2(0)$, \eqref{eq-II2} by $\omega'_1(0)$ and sum up both equations,  we get
\[
\left(\omega'_2(0) \overline{a} + \omega'_1(0) \overline{b}\right) \left(\overline{\alpha}\varphi_1+\varphi_2\right)\equiv 0\,,
\]
which implies (since $\overline{\alpha}\varphi_1+\varphi_2$ cannot be identically zero because otherwise, by \eqref{eq-II1}, we would have that $\varphi_3\equiv 0$ which is a contradiction, as was pointed out right before Lemma~\ref{lem-varphi})
\[
b=-\frac{\omega'_2(0)}{\omega'_1(0)}a=-pa\,.
\]
\par
Let us now replace those assertions in item ii) in Lemma~\ref{lem-recursive} in \eqref{eq-3}, and use \eqref{eq-identitycoefficients}, \eqref{eq-II1}, and \eqref{eq-II2} to get
 \[
\varphi_1 \overline{c\alpha}+\varphi_2 \overline{c} +\varphi_5\equiv 0\,.
\]
This gives, since $\varphi_5\not\equiv 0$ by Lemma~\ref{lem-varphi}, that $c\neq 0$.
\par
The substitution of  all the information obtained so far of the coefficients $a$, $b$, $c$, $d$, $k$, $l$, $m$, $n$, $\alpha$, and $\beta$  in equation \eqref{eq-I2} shows
\begin{equation}\label{eq-I2new}
ap\alpha\omega_1(\omega_2)^2 =\left(a\alpha(\omega_1)^2+c\alpha\omega_1-1\right)\omega_2 +(d\alpha+\beta)(\omega_1)^2+p\omega_1\,,
\end{equation}
while a re-arrangement of \eqref{eq-I3} gives
\begin{equation}\label{eq-I3new}
\left(1+ap\,\omega_1\right)(\omega_2)^2=\left(a(\omega_1)^2+c\omega_1\right)\omega_2+(p^2-cp)(\omega_1)^2\,.
\end{equation}
Now, multiply \eqref{eq-I2new} by $\left(1+ap\omega_1\right)$, multiply \eqref{eq-I3new} by $ ap\alpha\omega_1$, subtract the resulting equations and simplify to get
\[
\omega_2=\frac{p\omega_1+\left(a p^2+(p^2-cp)\alpha+\beta \right) (\omega_1)^2+ap\beta(\omega_1)^3}{1+\left(ap-c\alpha\right)\omega_1-a\alpha(\omega_1)^2}\,.
\]
We now substitute this expression for $\omega_2$ in \eqref{eq-afterprop-aux2},  multiply by $(1+\left(ap-c\alpha\right)\omega_1-a\alpha(\omega_1)^2)^2$ and re-arrange the terms to get an identity of the form
\[
a_1\omega_1+a_2(\omega_1)^2+a_3(\omega_1)^3+a_4(\omega_1)^4+a_5(\omega_1)^5+a_6(\omega_1)^6=0\,,
\]
where the coefficients $a_1, \ldots, a_6$ (which must be necessarily equal to zero) depend on the non-zero constants $a$, $c$, $p$, $\alpha$, and $\beta$.
\par
The coefficient $a_1=p\overline{\Phi_2(0)}-\overline{\Phi_1(0)}$ is directly equal to zero by \eqref{eq-P} and \eqref{eq-afterprop-aux1}. It is straightforward to check that the equation $a_2=0$ is automatically satisfied as well. A tedious but routine calculation shows that the expression $a_3=0$ is equivalent to the identity
\[
ap\,(p^2-1)=(c-2p)(\alpha p^2+\beta)\,,
\]
which gives, since $\alpha p^2+\beta=\beta(1-p^4)$,
\begin{equation}\label{eq-a}
ap=(2p-c)(1+p^2)\beta\,.
\end{equation}
Finally, from $a_5=0$ we obtain
\[
ap=2c\alpha +2 p(1+p^2)\beta\,.
\]
Hence, using also \eqref{eq-a} we have
\[
2\alpha+\beta(1+p^2)=0\,.
\]
It remains to replace the values of $p$, $\alpha$, and $\beta$ in terms of $\omega'_1(0)$ and $\omega'_2(0)$ in the previous equation to get the desired contradiction $\omega'_1(0)=\omega'_2(0)$.
\end{pf}
\section{The solution}\label{sec-solution}
\par\medskip
We now state the solution to the problem considered in this paper.
\begin{theorem}
Let $f$ be a locally univalent harmonic mapping in a simply connected domain $\Omega\subset \mathbb{C}$.
\par
\begin{itemize}
\item[i)] The only transformations that preserve local univalence and the harmonic Schwarzian derivative are pre-compositions with locally univalent affine harmonic mappings and anti-analytic rotations, in the cases when $f$ has non-constant dilatation.
\item[ii)] If the dilatation of $f$ is constant, so that $f=\alpha\overline{h}+h+\gamma$, where $h$ is an analytic function in $\Omega$, $|\alpha|\neq 1$, and $\gamma\in\mathbb{C}$ then, in addition to the transformations described in the previous item, we have that any other harmonic function $F$ of the form
    \[
    F=\beta\, \overline{T\circ h}+T\circ h+\lambda\,,
    \]
    where $|\beta|\neq 1$, $\delta\in\mathbb{C}$, and $T$ is a non-constant M\"{o}bius transformation as in \eqref{eq-Mobius} satisfies $S_H(f)=S_H(F)$.
\end{itemize}
\end{theorem}
\begin{pf}
The transformations considered preserve the Schwarzian derivative, as pointed out in Sections~\ref{ssec-properties} and \ref{ssec-equaldilat}.
\par
Conversely, if $f$ has non-constant dilatation and $F$ is another locally univalent harmonic function in $\Omega$ with $S_H(f)=S_H(F)$ we can compose both functions with a Riemann map from $\Omega$ onto the unit disk and apply a series of pre-compositions with (invertible) locally univalent affine harmonic mappings and anti-analytic rotations to get two new harmonic functions $\widetilde f$ and $\widetilde F$, say, now locally univalent and orientation-preserving in the unit disk and normalized as in \eqref{eq-normalizationh}, \eqref{eq-normalizationw}, and \eqref{eq-normalizationw'}. Then, by Theorem~\ref{thm-main}, the value of the derivative of the dilatations of these two functions at the origin must be equal. Hence, by Proposition~\ref{prop-equalfucn}, $\widetilde f=\widetilde F$. Undoing the invertible transformations used to produce $\widetilde f$ and $\widetilde F$, from $f$ and $F$, respectively, and using that the uniqueness principle holds for orientation-preserving harmonic mappings (see \cite[p. 8]{Dur-Harm}), we get the desired result.
\par
A direct application of Theorem~\ref{thm-constdilat1} ends the proof.
\end{pf}
%%%%%%%%%%%%%%%%%%%%%%%%%%%%%%%%%%%%%%%%%%%%%%%%%%%%%%%%%%%%%%%%%%%%%%

\end{document}